\documentclass[12pt]{article}

\usepackage{amsmath,amstext,amsfonts,amsthm,amssymb,eucal,diagrams}

\newcommand{\bzero}{{\bar{0}}}
\newcommand{\bone}{{\bar{1}}}

\newcommand{\be}{\begin{equation}}
\newcommand{\ee}{\end{equation}}

\renewcommand{\mod}{\operatorname{mod}}
\newcommand{\Com}{\operatorname{Com}}

\newcommand{\Tor}{\operatorname{Tor}}

\newcommand{\OO}{{\cal O}}

\newcommand{\coker}{\operatorname{coker}}

\newcommand{\DD}{{\cal D}}

\newcommand{\RR}{{\cal R}}
\newcommand{\KK}{{\cal K}}

\newcommand{\G}{{\mathbb G}}

\newcommand{\hra}{\hookrightarrow}

\newcommand{\HHom}{{\cal H}om}

\newcommand{\Coh}{\operatorname{Coh}}
\newcommand{\GG}{{\cal G}}
\newcommand{\CC}{{\cal C}}

\newcommand{\Spec}{\operatorname{Spec}}
\newcommand{\Sing}{\operatorname{Sing}}

\newcommand{\perf}{\operatorname{Perf}}

\newcommand{\ga}{\gamma}
\newcommand{\de}{\delta}

\renewcommand{\ker}{\operatorname{ker}}

\newcommand{\D}{{\cal D}}

\numberwithin{equation}{section}

\newtheorem{theor}{Theorem}[section]
\newtheorem{thm}[theor]{Theorem}
\newtheorem{lem}[theor]{Lemma}

\newtheorem{prop}[theor]{Proposition}

\newtheorem{cor}[theor]{Corollary}
{  \theoremstyle{definition}
           \newtheorem{defi}[theor]{Definition}
           \newtheorem{rem}[theor]{Remark}

}

\newcommand{\Pf}{\noindent {\it Proof}}
\newcommand{\id}{\operatorname{id}}

\newcommand{\ov}{\overline}

\newcommand{\FF}{{\cal F}}

\newcommand{\PP}{{\cal P}}

\newcommand{\Hom}{\operatorname{Hom}}

\newcommand{\Ext}{\operatorname{Ext}}

\renewcommand{\a}{\alpha}

\newcommand{\la}{\lambda}

\newcommand{\Z}{{\mathbb Z}}

\newcommand{\Ga}{\Gamma}

\newcommand{\wt}{\widetilde}
\newcommand{\ot}{\otimes}

\newcommand{\sub}{\subset}
\newcommand{\com}{\operatorname{com}}
\newcommand{\ed}{\qed\vspace{3mm}}

\newcommand{\Qcoh}{\operatorname{Qcoh}}

\newcommand{\Sg}{\operatorname{Sg}}

\newcommand{\MF}{{\operatorname{MF}}}
\newcommand{\HMF}{\operatorname{HMF}}
\newcommand{\DMF}{\operatorname{DMF}}
\newcommand{\mf}{matrix factorization}

\newcommand{\Homb}{\mathcal{H}om} 
\newcommand{\Per}{\operatorname{Per}}
\newcommand{\fC}{{\mathfrak C}}
\newcommand{\LHZ}{\operatorname{LHZ}}

\newcommand{\Lfr}{\operatorname{Lfr}}
\newcommand{\EExt}{\mathcal{E}xt}

\title{Matrix factorizations and singularity categories for stacks}
\author{Alexander Polishchuk \and Arkady Vaintrob}
\date{}

\begin{document}
\maketitle
\begin{abstract}
We study \mf s of a section $W$ of a line bundle on an algebraic stack.
 We relate the corresponding derived category (the category of D-branes of type B in
 the Landau-Ginzburg model with potential $W$)
 with the singularity category of the zero locus of 
 $W$ generalizing a theorem of Orlov.
 We use this result to construct push-forward functors 
 for \mf s with relatively proper support.
\end{abstract}

\section*{Introduction}

Matrix factorizations arose in the work of Eisenbud \cite{Eisen}
in connection with the study of maximal Cohen-Macaulay modules.
Since then they became a standard tool and an object of study in commutative algebra 
(see  e.g.\ \cite{BGS, Yoshino}).
Recall that a \mf\ of an element $W$ of a commutative ring $R$ is a $\Z/2$-graded
finitely generated projective $R$-module $E=E_0\oplus E_1$
equipped with an odd endomorphism $\de:E\to E$ such that $\de^2=W\cdot\id$.
Matrix factorizations of 
a fixed element $W$ (called a \emph{potential})
form a triangulated category $\HMF(W)$ 
with morphisms defined as ``chain maps'' up to homotopy.

Following the suggestion of Kontsevich, \mf s were used
by physicists to describe D-branes of type B in Landau-Ginzburg models
(see \cite{KL1, KL2}).
They found applications in various approaches to mirror symmetry 
and   in the study of sigma model/Landau-Ginzburg correspondence
(see \cite{KKP, Ef, Seidel, BHLW, HHP}).
Mathematical foundations of this circle of ideas were laid down in a series of papers by Orlov
        \cite{Orlov,Orlov-graded,O-compl}. 

The fundamental result of Orlov \cite[Thm.\ 3.9]{Orlov} states that the triangulated
category $\HMF(W)$ is equivalent to the so-called {\it singularity category}
$D_{\Sg}(X_0)$ of the affine hypersurface $X_0=\Spec(R/(W))\sub X=\Spec(R)$ 
(assuming that $X$ is smooth and $W$  is not a zero divisor).
Here $D_{\Sg}(X_0)$ is defined as the quotient of the bounded derived category $D^b(X_0)$
of coherent sheaves on $X_0$ by the triangulated subcategory $\Per(X_0)$ of perfect complexes.
One should think of a triangulated category $D_{\Sg}(X_0)$ as a ``measure'' of singularity of $X_0$
(in the case when $X_0$ is smooth one has $D_{\Sg}(X_0)=0$).

Kontsevich (see \cite{Kon-lec}) suggested to view categories of \mf s as examples of noncommutative
spaces. In this context it is natural to consider analogs of these categories for a potential $W$
on a non-necessarily affine scheme $X$. 
In this case simply using the homotopy category $\HMF(W)$ does not give the right notion
because one has to take into account the presence of higher cohomology of coherent sheaves on 
$X$. One way to do this is to use a dg-version of the functor $R\Hom$
(see \cite[Sec.\ 3.2]{KKP}, \cite{Segal-GIT}). Another (proposed by Orlov, see \cite[Sec.\ 3.2]{KKP}) 
is to replace the homotopy category of \mf s by an appropriate localization.
This is the path we follow in this paper. We define the derived category of \mf s $\DMF(X,W)$
as the quotient of the homotopy category $\HMF(X,W)$ by the triangulated subcategory formed
by \mf s that are locally contractible. 
This definition of $\DMF(X,W)$ is different from the one proposed by Orlov (see Remark \ref{Orlov-def})
but is equivalent to it under appropriate assumptions on $X$.

In examples coming from physics it is often necessary
consider categories of \mf s on orbifolds. 
To include this case we consider a more general situation where $X$ is an algebraic stack
and $W$ is a section of a line bundle over $X$. Our main result  (see Theorem
\ref{equiv-thm}) is a generalization  of Orlov's theorem to this
case. Namely, we show that for a potential $W$ on a smooth algebraic stack
$X$ satisfying some technical assumptions the derived category $\DMF(X,W)$ of
\mf s is equivalent to the singularity category of the zero locus of $W$. 

The second goal of the paper is to define and study push-forward functors for \mf s
as a preparation to developing an analog of the theory of Fourier-Mukai transforms.
The naive definition of the push-forward functor with respect to a smooth affine morphism
leads to \mf s of possibly infinite rank (we call them {\it quasi-\mf s}).
On the other hand, using the equivalence with the singularity category and the notion of
support for \mf s which we develop in Section \ref{supports-sec}, we define a
push-forward functor for derived category of \mf s with relatively proper support. 
To check that these two types of functors are compatible
we prove some partial analogs of Orlov's equivalence for quasi-\mf s.

In the paper \cite{PV-CohFT}
we will use the results of this paper to provide an algebraic analog
of the theory developed by Fan, Jarvis and Ruan in \cite{FJR, FJR2}. For each 
quasihomogeneous polynomial $W$ with an isolated singularity
they construct a    cohomological field theory which is a Landau-Ginzburg counterpart of the 
topological sigma model producing Gromov-Witten invariants.
The main technical ingredient of the Fan-Jarvis-Ruan theory 
is a collection of cohomology classes on the moduli space 
of stable curves with additional data (called $W$-curves).
The construction given in \cite{FJR2} uses sophisticated analytic tools.
We construct in \cite{PV-CohFT} similar classes using categories 
of \mf s and functors between them introduced in this paper.
The starting point of our approach is the identification 
made in our previous paper \cite{PV-mf} of the orbifold Milnor ring of $W$ 
with its residue paring (which is equal to the state space of the
Fan-Jarvis-Ruan theory) with the Hochschild homology of the dg-category 
of equivariant \mf s of $W$ equipped with the canonical metric.

Now let us describe the structure of the paper.
In Section \ref{1st-sec} we give basic definitions and constructions for \mf s on stacks.

The new feature here is that we allow the potential $W$ to be a section of a
non-trivial line bundle $L$ on $X$.
In this setup we have a natural definition of $\Z$-graded dg-categories of \mf s, 
making our \mf s similar to graded \mf s considered in \cite{Orlov-graded, KKT}.
The $\Z/2$-graded dg-categories of \mf s (considered originally in \cite{Orlov}) can be considered
only in the case when $L$ is trivial.
In Section \ref{equiv-sec} we show how our formalism can be used to 
work with equivariant \mf s. 
Our main result, Theorem \ref{equiv-thm}, generalizing Orlov's equivalence 
\cite[Thm.\ 3.9]{Orlov} between singularity categories and derived
categories of matrix factorizations to the non-affine case
is proved in Section \ref{sing-sec}. 
As in the affine situation, there is a natural
functor from the homotopy category $\HMF(W)$ to the singularity category $D_{\Sg}(X_0)$
of the zero locus of $W$. The main difficulty in the non-affine case is to prove
the surjectivity of the induced functor from the derived category $\DMF(W)$ to $D_{\Sg}(X_0)$.
This is achieved by combining Orlov's description of morphisms in $D_{\Sg}(X_0)$ (see
\cite[Prop.\ 1.21]{Orlov}) with the standard $2$-periodic resolution associated with a \mf.
In Section \ref{quasi-sec} we prove some partial analogs
of this equivalence for quasi-\mf s  
--- factorizations of possibly infinite rank which will be used in the study of push-forward functors.
In Section \ref{supports-sec} we introduce
the notion of support for \mf\ s. 
By definition, the support of $(E,\de)$ is the subset of
the zero locus of $W$ consisting of points $x$ such that the $2$-periodic complex
defined by the pull-back of $(E,\de)$ to $x$ has nonzero cohomology.
We prove that this corresponds to the natural notion of support for objects of the singularity
category. More precisely, we have to work with the idempotent completion of this category (following
\cite{O-compl}).
Finally, in Section \ref{push-forward-sec}
we study push-forward functors for categories of \mf s.
We define the push-forward of \mf s with relatively proper support
for a representable morphism of smooth stacks (satisfying additional technical assumption),
taking values in the idempotent completion of the derived category
of \mf s on the base. We prove that for a smooth affine morphism with
geometrically integral fibers this push-forward agrees with the obvious notion of a push-forward
for quasi-matrix factorizations.

\

{\it Acknowledgments}. 
We would like to thank Maxim Kontsevich, Dmitry Orlov and Tony Pantev for helpful discussions. 
Part of this work was done during our stay at the IHES and we are grateful to it
for hospitality and for stimulating atmosphere.
The first author was partially supported by the NSF grant DMS-1001364.

\

{\it Notations and conventions}.
We work with schemes and stacks over a ground field $k$.
All stacks we are dealing with are assumed to be algebraic, Noetherian and semi-separated.
For such a stack $X$ we denote by $\Coh(X)$, (resp., $\Qcoh(X)$; resp., $D^b(X)$) the category of coherent sheaves (resp., quasicoherent sheaves; resp., bounded derived
category of coherent sheaves) on $X$. By \cite[Cor.\ 2.11]{AB}, $D^b(X)$ 
is equivalent to the full subcategory
of the bounded derived category of $\OO_X$-modules consisting of complexes with coherent
cohomology.
We call a stack $X$ Gorenstein if there exists a presentation $U\to X$, where $U$ is
a Gorenstein scheme.
This implies that the structure sheaf $\OO_X$ is a dualizing object in $D^b(X)$.
By vector bundles we mean locally free sheaves of finite rank.

\section{Matrix factorizations of a section of a line bundle}
\label{1st-sec}

Let $X$ be an algebraic stack, 
$L$ a line bundle on $X$, and $W\in H^0(X,L)$ a section (called a {\it potential}). 

\begin{defi}
A {\it \mf\ } $\bar{E}=(E_\bullet,\de_\bullet)$ of $W$ on $X$ consists of a pair of vector bundles 
(i.e., locally free sheaves of finite rank) $E_0$, $E_1$
on $X$ together with homomorphisms 
$$\de_1:E_1\to E_0 \ \text{ and } \ \de_0:E_0\to E_1\ot L,$$ 
such that
$\de_0\de_1=W\cdot \id$ and $\de_1\de_0=W\cdot\id$.

We will often assume that the potential $W$ is not a zero divisor, i.e., the morphism
$W:\OO_X\to L$ is injective. 
\end{defi}

It is convenient to introduce formal expressions $V\ot L^{n/2}$, where $V$ is a vector bundle $X$,
which we will call {\it half-twisted bundles}.
They have a natural tensor product 
$$(V_1\ot L^{n_1/2})\ot (V_2\ot L^{n_2/2})=V_1\ot V_2\ot L^{(n_1+n_2)/2}.$$
The space of morphisms between half-twisted bundles $V_1\ot L^{n_1/2}$ and
$V_2\ot L^{n_2/2}$
is defined only
when $n_2-n_1$ is even and given by
$$\Hom(V_1\ot L^{n_1/2}, V_2\ot L^{n_2/2})=\Hom_X(V_1,V_2\ot L^{(n_2-n_1)/2}).$$ 
We will also use $\Z/2$-graded half-twisted bundles
$V=\left(V_\bzero\ot L^{m/2}\right)\oplus \left(V_\bone\ot L^{n/2}\right)$ and we
define their ``half-twists" by
$$V(L^{1/2}):=V_\bzero\ot L^{m/2}\oplus V_\bone\ot L^{(n+1)/2} \ \text{ and } 
V(L^{-1/2}):=V_\bzero\ot L^{m/2}\oplus V_\bone\ot L^{(n-1)/2}.$$

With a \mf\ $\bar{E}=(E_\bullet,\de_\bullet)$ we associate a $\Z/2$-graded half-twisted bundle
$$E(L^{1/2})=E_0\oplus (E_1\ot L^{1/2}).$$ 
The differential
$\de$ can be viewed as an odd morphism 
$$\de:E(L^{1/2})\to E(L^{1/2})\ot L^{1/2}$$ such that $\de^2=W$.

\begin{defi} We define the dg-category $\MF(X,W)$ of  \mf s of $W$ as follows.
For \mf s $\bar{E}$ and $\bar{F}$ the set of morphisms
$\Homb_\MF(\bar{E},\bar{F})$ is  a $\Z$-graded complex 
$$\Homb_\MF(\bar{E},\bar{F})^i=\Hom_{i\mod 2}(E(L^{1/2}),F(L^{1/2})\ot L^{i/2}),$$
where $\Hom_\bzero$ (resp., $\Hom_\bone$) denote the even (resp., odd) morphisms
between the $\Z/2$-graded half-twisted bundles.
Explicitly, 
$$\Homb_\MF(\bar{E},\bar{F})^{2n}=\Hom(E_0,F_0\ot L^n)\oplus\Hom(E_1,F_1\ot L^n),$$
$$\Homb_\MF(\bar{E},\bar{F})^{2n+1}=\Hom(E_0,F_1\ot L^{n+1})\oplus\Hom(E_1,F_0\ot L^n).$$
The differential on $\Homb_\MF(\bar{E},\bar{F})$ is given by
\begin{equation}\label{diff-mf-eq}
d f = \de_F\circ f - (-1)^{|f|} f \circ \de_E~.
\end{equation}
\end{defi}

We denote by $\HMF(X,W)=H^0\MF(X,W)$ the corresponding homotopy category. In other words,
the objects of $\HMF(X,W)$ are \mf s of $W$, morphisms are the spaces
$H^0 \Homb_\MF(\bar{E},\bar{F})$.
We will usually omit $X$ from the notation.
As in the affine case considered in \cite{Orlov}, the category $\HMF(W)$ has a triangulated structure.

\begin{defi}\label{triangle-def} 
We define
the translation functor on $\HMF(W)$ by
$$\bar{E}[1]=(E[1],\de[1]), \ \text{ where } E[1]_0=E_1\ot L, E[1]_1=E_0, \de[1]_i=-\de_{i+1}.$$
The mapping cone of a closed morphism of \mf s $f:\bar{E}\to\bar{F}$
is defined as $C(f)=F\oplus E[1]$ with the differential given
by the same formula as for the category of complexes. We have canonical closed morphisms of \mf s
$\bar{F}\to C(f)$ and $C(f)\to\bar{E}[1]$. We define the class of exact triangles
as those isomorphic to some triangle of the form 
\begin{equation}\label{cone-triangle}
\bar{E}\rTo{f}\bar{F}\to C(f)\to\bar{E}[1].
\end{equation}
\end{defi}

The standard proof that the homotopy category of complexes is triangulated can be adapted
to show that we obtain a triangulated structure on $\HMF(W)$ in this way.

Note that for $n\in\Z$ one has
$$E[n](L^{1/2})_i=E(L^{1/2})_{i+n}\ot L^{n/2}$$
with $\de[n]_i=(-1)^n\de_{i+n}$.

Sometimes it is convenient to work with \mf s of possibly infinite rank. 

\begin{defi} A {\it quasi-\mf\ } $\bar{E}=(E_0\oplus E_1,\de)$ of $W\in H^0(X,L)$
consists of a pair of locally free
sheaves (not necessarily of finite rank) $E_0$ and $E_1$ equipped with the differentials
$\de_0,\de_1$ as above. As before, quasi-\mf s form
a dg-category $\MF^{\infty}(X,W)$, and the corresponding homotopy category 
$\HMF^{\infty}(X,W)$ is triangulated.
\end{defi}

Matrix factorizations can be used to produce infinite ($2$-periodic up to a twist) complexes
of vector bundles on the zero locus $X_0=W^{-1}(0)$. Namely, with
a quasi-\mf\ $\bar{E}=(E_\bullet,\de_\bullet)$ we associate a $\Z$-graded complex
$\com(\bar{E})$ of vector bundles on $X_0$
\begin{equation}\label{mf-Z-gr-com-eq}
\ldots\to (E_0\ot L^{-1})|_{X_0} \rTo{\ov{\de}_0} E_1|_{X_0}\rTo{\ov{\de}_1} E_0|_{X_0}\rTo{\ov{\de}_0}
(E_1\ot L)|_{X_0}\to\ldots
\end{equation}
where $\ov{\de}_i$ is induced by $\de_i$, and $E_0|_{X_0}$ is placed in degree $0$.
Note that we can also present $\com(\bar{E})$ using half-twisted bundles:
$$\com(\bar{E})^n=E(L^{1/2})_n\ot L^{n/2}|_{X_0}.$$
This construction extends to a dg-functor
$$\com:\MF^\infty(W)\to\Com(X_0)$$
that induces an exact functor 
$\com:\HMF^\infty(W)\to H^0\Com(X_0)$
between the corresponding homotopy categories.

\begin{lem}\label{2-per-res-lem} If $W$ is not a zero divisor,
then for any quasi-\mf\ $\bar{E}$ of $W$ the complex $\com(\bar{E})$ is exact.
\end{lem}

\Pf . For a coherent sheaf $\FF$ on $X$ let us denote by $W\FF\sub \FF$ the image of the map
$\FF\ot L^{-1}\to \FF$ induced by $W$. The kernel of the map
$\ov{\de}_1:E_1/WE_1\to E_0/WE_0$
can be identified with $\de_1^{-1}(WE_0)/WE_1$.
But $WE_0=\de_1\de_0E_0$ and $\de_1$ is injective, so we obtain
$\de_1^{-1}(WE_0)=\de_0E_0$. Now $\de_0$ induces an isomorphism
$$E_0/\de_1 E_1\simeq \de_0(E_0\ot L^{-1})/WE_1\simeq \ker(\ov{\de}_1).$$
But 
$$E_0/\de_1 E_1\simeq\coker(\de_1)\simeq\coker(\ov{\de}_1:E_1/WE_1\to E_0/WE_0),$$
which finishes the proof.
\ed

For a morphism of stacks $f:X'\to X$, a line bundle $L$ over $X$
and a section $W\in H^0(X,L)$ we have natural pull-back functors of \mf s: 
a dg-functor
$$f^*:\MF(X,W)\to\MF(X',f^*W),$$
where $f^*W$ is the induced section of $f^*L$ on $X'$, and the induced exact functor
$$f^*:\HMF(X,W)\to\HMF(X,f^*W).$$

\section{Equivariant \mf s}
\label{equiv-sec}

Here we define equivariant \mf s and show that they can be viewed as a particular case of
the construction of the previous section.

Let $\Ga$ be an affine algebraic
group acting on a stack $X$. Let $W$ be a regular function on $X$, semi-invariant with respect to $\Ga$,
i.e., we have a character $\chi:\Ga\to\G_m$ such that 
$$W(\ga\cdot x)=\chi(\ga)W(x)$$
for $\ga\in \Ga$, $x\in X$. 

\begin{defi}
A {\it $\Ga$-equivariant \mf\ of $W$ with respect to the character } $\chi$
is a pair of $\Ga$-equivariant vector bundles $(E_0, E_1)$ on $X$ together with
$\Ga$-invariant homomorphisms
$$\de_1:E_1\to E_0\ \text{ and } \de_0:E_0\to E_1\ot\chi,$$
such that $\de_0\de_1=W\cdot \id$ and
$\de_1\de_0=W\cdot \id$. 

$\Ga$-equivariant \mf s form a dg-category $\MF_{\Ga,\chi}(X,W)$ with
morphisms given by
$$\Homb_{\MF_\Ga}(\bar{E},\bar{F})^{2n}=\Hom(E_0,F_0\ot \chi^n)^\Ga\oplus\Hom(E_1,F_1\ot \chi^n)^\Ga,$$
$$\Homb_{\MF_\Ga}(\bar{E},\bar{F})^{2n+1}=\Hom(E_0,F_1\ot \chi^{n+1})^\Ga\oplus\Hom(E_1,F_0\ot \chi^n)^\Ga$$
and the differential defined by \eqref{diff-mf-eq}.
The corresponding homotopy category $\HMF_{\Ga,\chi}(X,W)$ has a triangulated structure defined
as in Definition \eqref{triangle-def}. We will usually omit $X$ and $\chi$ from the notation.
\end{defi}

Note that closed morphisms of degree zero between
$\Ga$-equivariant \mf s
$\bar{E}$ and $\bar{F}$ are given by pairs of
$\Ga$-invariant morphisms $(E_0\to F_0, E_1\to F_1)$ commuting with
$\de$. 

Equivariant \mf s of $W$ can be described 
as \mf s on the stack $Y=X/\Ga$ with respect to the line bundle
$L_{\chi^{-1}}$ which corresponds to the $\G_m$-torsor defined as
the push-out of the $\Ga$-torsor $U\to X$ under $\chi^{-1}$.
Note that $W$ descends to a global section $\ov{W}\in H^0(Y,L_{\chi^{-1}})$.
Viewing vector bundles on $Y$ as $\Ga$-equivariant vector bundles on $X$ we immediately
obtain the following result.

\begin{prop}\label{equiv-mf-prop}
The dg-categories (resp., triangulated categories) $\MF_{\Ga,\chi}(X,W)$ and
$\MF(Y,\ov{W})$ (resp., $\HMF_{\Ga,\chi}(X,W)$ and $\HMF(Y,\ov{W})$) are equivalent.
\end{prop}

Conversely, if $W$ is a section of a line bundle $L$ on a stack $X$, then \mf s of $W$ can
be viewed as $\G_m$-equivariant \mf s on the $\G_m$-torsor $\wt{X}\to X$
associated with $L^{-1}$. Indeed, $W$ can be viewed as a function
$\wt{W}$ on $\wt{X}$, such that $\wt{W}(\la\cdot \wt{x})=\la\cdot\wt{W}(\wt{x})$ for $\la\in\G_m$. In other words, $\wt{W}$ is semi-invariant with respect to $\G_m$ and the character 
$\id:\G_m\to\G_m$.
Proposition \ref{equiv-mf-prop} applied to the action of $\G_m$ on $\wt{X}$ implies
that the dg-categories 
$\MF(X,W)$ and $\MF_{\G_m,\id}(\wt{X},\wt{W})$ are equivalent.

Let $X$ (resp., $X'$) be a stack with an action of a group $\Ga$
(resp., $\Ga'$). Assume that we have a homomorphism $\pi:\Ga'\to \Ga$ and a $\Ga'$-equivariant
morphism $f:X'\to X$. Given a character $\chi:\Ga\to\G_m$ and a function
$W$ on $X$ such that $W(\ga\cdot x)=\chi(\ga)W(x)$, let $f^*W$ be the pull-back of $W$ to $X'$.
Then $f^*W$ is semi-invariant with respect to the character $\chi'=\chi\circ\pi$ of $\Ga'$.
In this situation we have a natural pull-back 
dg-functor
\be\label{pull-back-ex} 
\MF_{\Ga,\chi}(X,W)\to\MF_{\Ga',\chi'}(X',f^*W).
\end{equation}

\section{Connection with categories of singularities}
\label{sing-sec}

In this section we define the derived category $\DMF(X,W)$ of \mf s of a potential $W\in H^0(X,L)$
on a stack $X$ as a suitable localization of $\HMF(X,W)$. Our main result is Theorem \ref{equiv-thm}
establishing equivalence of $\DMF(X,W)$ with the singularity category of $X_0$, the zero
locus of $W$. This is a generalization of \cite[Thm.\ 3.9]{Orlov}.

\begin{defi}\label{stack-def} 
Let $X$ be an algebraic stack. 

\noindent
We say that $X$ 

(i) has the {\it resolution property (RP)} if
for every coherent sheaf $\FF$ on $X$ there exists a vector bundle $V$ on $X$
and a surjection $V\to\FF$.

(ii) has {\it finite cohomological dimension (FCD)} if 
there exists an integer $N$ such that for every quasicoherent sheaf $\FF$ on $X$ one
has $H^i(X,\FF)=0$ for $i>N$. We call the minimal $N$ with this property the {\it cohomological
dimension} of $X$.

(iii) is an {\it FCDRP-stack} if it has both properties (i) and (ii).

(iv) is a {\it nice quotient stack} if $X=U/\Ga$,
where $U$ is a Noetherian scheme and $\Ga$ is a reductive linear
algebraic group, such that $U$ has an ample family of $\Ga$-equivariant line bundles.
\end{defi}

Note that by \cite{Thom}, if $X$ is a nice quotient stack then it has the resolution property.
Such a stack also has finite cohomological dimension, as the following lemma shows.

\begin{lem}\label{PF-lem} 
Let $X=U/\Ga$ be a quotient stack, where $U$ is a scheme and $\Ga$ is a reductive group,
and let $\pi:U\to X$ be the natural projection.

(i) For any quasicoherent sheaves $\FF$ and $\GG$ on $X$ one has
natural isomorphisms
$$\Ext^i_X(\FF,\GG)\simeq\Ext^i_U(\pi^*\FF,\pi^*\GG)^\Ga.$$

(ii) Assume that $U$ can be covered with $N$ affine open subsets.
If $\PP$ is a locally free sheaf on $X$ and $\FF$ is a quasicoherent sheaf on $X$ then
$$\Ext^i_X(\PP,\FF)=0 \ \text{ for } i>N.$$
In particular, $X$ has cohomological dimension $\le N$.
\end{lem}

\Pf . (i) For $i=0$ this follows from the definition of the morphisms in the category of $\Ga$-equivariant
sheaves on $U$. The general case follows from the exactness of the functor of $\Ga$-invariants.

(ii) By part (i), it suffices to prove a similar fact on the scheme $U$, which follows from the proof of
\cite[Lem.\ 1.12]{Orlov}.
\ed
 
\begin{rem}
By \cite[Prop.\ 5.1]{Kresch} if $X$ is a separated DM stack of finite type over $k$ 
that has quasi-projective coarse moduli space then $X$ has RP if and only it is a
quotient stack. Moreover, in this case $X$ is a nice quotient stack.
For example, this is true for every smooth separated DM stack of finite type over a field of characteristic
zero with quasi-projective coarse moduli space (see \cite[Thm.\ 4.4]{Kresch}).
\end{rem}

The condition of finite cohomological dimension is stable under passing to open and closed substacks.
More generally, we have the following simple observation.

\begin{lem}\label{FCD-stack-lem}
Let $f:X\to Y$ be a representable morphism of Noetherian stacks. 
If $Y$ has finite cohomological dimension
then so does $X$.
\end{lem}

\Pf . It suffices to prove that $Rf_*$ has finite cohomological dimension. But this can be checked
by replacing $Y$ with its presentation, so we can assume $Y$ to be a scheme. Now
the assertion follows from the assumption that our
stacks are Noetherian.
\ed
 
Recall (see \cite{SGA6I}) that the full subcategory $\Per(X)\in D^b(X)$ of {\it perfect complexes} 
consists of objects, locally isomorphic in the derived category
to a bounded complex of vector bundles.

\begin{lem}\label{per-lem} Let $X$ be a stack with the resolution property.

(i) For every complex $C^\bullet\in D^b(X)$ there exists a bounded above complex of vector bundles $P^\bullet$ and a quasi-isomorphism $P^\bullet\to C^\bullet$.

(ii) Any object of $\Per(X)$ is isomorphic in $D^b(X)$ 
 to a bounded complex of vector bundles.

(iii)  For every object
$C^\bullet\in\Per(X)$ there exists a bounded complex of vector bundles
$P^{\bullet}$ and a quasi-isomorphism $P^\bullet\to C^\bullet$.
\end{lem}

\Pf . (i) and (ii) follow from the resolution property
as in \cite[Lem.\ 2.2.8]{SGA6II} using \cite[Prop.\ 1.2]{SGA6II}.
As is explained in \cite[Lem.\ 1.6]{Orlov}, (iii) follows from (i) (one can also apply 
of \cite[Prop.\ 2.2]{SGA6I}).
\ed

\begin{defi}
Let $X$ be a stack.
The quotient category
$$D_{\Sg}(X):=D^b(X)/\Per(X)$$ 
is called the {\it singularity category of } $X$.
We denote by $\ov{D_{\Sg}(X)}$ its idempotent completion.
\end{defi}

The fact that being a perfect complex is a local property immediately implies the following result.

\begin{lem}\label{quasi-isom-lem} 
An object $A\in D_{\Sg}(X)$ (resp., a morphism $f:A\to B$ in $D_{\Sg}(X)$) 
is zero (resp., an isomorphism) if and only if
it is zero (resp., an isomorphism) locally, i.e., for some open covering $(U_i)$ of $X$ in flat topology
the restrictions $A|_{U_i}$) are zero (resp., the restrictions
$f|_{U_i}$ are isomorphisms in $D_{\Sg}(U_i)$).
\end{lem}

The singularity category $D_{\Sg}(X)$ admits a natural quasicoherent analog. 
Namely, we consider the full subcategory $\Lfr(X)\sub D^b(\Qcoh(X))$ of objects that can be represented
by bounded complexes of locally free sheaves, and consider the quotient
\begin{equation}\label{quasi-sing-eq}
D'_{\Sg}(X):= D^b(\Qcoh(X))/\Lfr(X).
\end{equation}
Note that if $X$ has a resolution property then by Lemma \ref{per-lem}(ii), 
we have $\Lfr(X)\cap D^b(X)=\Per(X)$.

\begin{prop}\label{quasicoh-emb-lem} Let $X$ be a nice quotient stack 
(see Definition \ref{stack-def}).
Then the natural functor $D_{\Sg}(X)\to D'_{\Sg}(X)$ is fully faithful.
\end{prop}

\Pf . We can repeat the argument of \cite[Prop.\ 1.13]{Orlov} using Lemma \ref{PF-lem}(ii) and Lemma
\ref{per-lem}.
\ed

\begin{defi} A coherent sheaf $\FF$ on a Gorenstein stack $X$ of finite Krull dimension
is called maximal Cohen-Macaulay (MCM) if
$R\Homb(\FF,\OO_X)$ is a sheaf.
\end{defi} 

\begin{lem}\label{FCD-lem} 
Let $X$ be a Gorenstein stack of finite cohomological dimension and finite Krull
dimension. Then there exists an integer $N$ such that
for any coherent sheaf $\FF$ and any locally free sheaf $\PP$ on $X$ one has
$\Ext^i(\FF,\PP)=0$ for $i>N$.
\end{lem}

\Pf . The proof is similar to that of \cite[Lem.\ 1.18]{Orlov}.
Namely, first, one checks that $\EExt^i(\FF,\PP)=0$ for $i>n$, where $n$ is the Krull dimension of $X$,
using the fact that $X$ is Gorenstein. 
Next, since the sheaves $\EExt^i(\FF,\PP)$ are quasicoherent, the result follows from the local to global spectral sequence and the assumption that $X$ has finite cohomological dimension.
\ed

\begin{lem}\label{MCM-lem} Let $X$ be a Gorenstein stack of finite Krull dimension
with the resolution property.

(i) A coherent sheaf $\FF$ on $X$
is MCM if and only if it admits an (infinite) right resolution by vector bundles.

(ii) A MCM sheaf, which is perfect as a complex, is locally free.

(iii) Every object of $D_{\Sg}(X)$ is isomorphic
to a MCM sheaf.

(iv) Assume in addition that $X$ has finite cohomological dimension $N$.
Let $\FF$ be a MCM sheaf, $\GG$ a quasicoherent sheaf on $X$. 
Then for any morphism $f:\FF\to\GG$ in $D'_{\Sg}(X)$,
any integer $n\ge N$ and any exact sequence of coherent sheaves on $X$
\begin{equation}\label{n-ex-seq1}
0\to\GG'\to G_n\to\ldots\to G_1\to\GG\to 0,
\end{equation}
where $G_i$ are locally free and $\GG'$ is a quasicoherent sheaf on $X$,
there exists an exact sequence of coherent sheaves
\begin{equation}\label{n-ex-seq2}
0\to \FF'\to F_n\to\ldots\to F_1\to \FF\to 0,\ \ 
\end{equation}
where $F_i$ are vector bundles on $X$, and a morphism $f':\FF'\to\GG'$ of quasicoherent sheaves
making the following diagram in $D'_{\Sg}(X)$ is commutative
\begin{diagram}
\FF &\rTo{}& \FF'[n]\\
\dTo{f}&&\dTo{f'[n]}\\
\GG&\rTo{}&\GG'[n]
\end{diagram}
(here the horizontal arrows are isomorphisms in $D'_{\Sg}(X)$ induced by the exact sequences
\eqref{n-ex-seq1}, \eqref{n-ex-seq2}).
\end{lem}

\Pf . (i) See  \cite[Lem.\ 1.19]{Orlov}. 

(ii) This is proved in the same way as \cite[Lem.\ 1.20]{Orlov}. If $\FF$ is such a sheaf then
by Lemma \ref{per-lem}(iii), there exists a bounded complex of vector bundles $P^\bullet$
and a quasi-isomorphism $P^\bullet\to \FF^\vee=\HHom(\FF,\OO_X)$ (since
$\FF^\vee$ is perfect). Dualizing we obtain a bounded right resolution for $\FF$, so $\FF$
is locally free.

(iii) One can repeat the proof of \cite[Prop.\ 1.23]{Orlov} using Lemma \ref{per-lem}(i) instead
of \cite[Lem.\ 1.4]{Orlov}.

(iv) We can view $f$ as a morphism $\FF\to \GG'[n]$ in $D'_{\Sg}(X)$. 
We can apply the same argument as in \cite[Prop.\ 1.21]{Orlov}
to prove that $f$ can be realized by a morphism $\FF\to\GG'[n]$ in $D^b(\Qcoh(X))$.
Note that because $\GG'$ is only quasicoherent, and $D_{\Sg}(X)$ is replaced by $D'_{\Sg}(X)$,
we have to use Lemma \ref{FCD-lem} at the relevant place in this argument.
Finally, as is well known, a morphism $\FF\to\GG'[n]$ in $D^b(\Qcoh(X))$
can be realized by a morphism $H_{n}(F_\bullet)\to\GG'$
for appropriate resolution of $F_\bullet$ of $\FF$. Furthermore, since $X$ has the resolution
property, we can assume that $F_i$ are vector bundles.
\ed

Now let $X$ be a stack, and let
$W\in H^0(X,L)$ be a potential, where $L$ is a line bundle on $X$.
Assume that $W$ is not a zero divisor.
Let $X_0=W^{-1}(0)$ be the zero locus of $W$. As in \cite{Orlov}, we consider the natural functor 
\begin{equation}\label{coker-functor-eq}
\fC:\HMF(X,W)\to D_{\Sg}(X_0).
\end{equation}
that associates with a \mf\ $(E_{\bullet},\de)$ the cokernel of $\de_1:E_1\to E_0$. 

\begin{lem}\label{exact-lem}
The functor $\fC$ is exact.
\end{lem}

\Pf . For a \mf\ 
$$\bar{F}=(F_0\oplus F_1,\de)$$
consider
$$\FF=\fC(\bar{F})=\coker(F_1\rTo{\de_1} F_0)=\coker(F_1|_{X_0}\rTo{\ov{\de}_1} F_0|_{X_0})\ \ \text{ and}$$
$$\FF'=\fC(\bar{F}[1])=\coker(F_0\rTo{\de_0} F_1\ot L)=\coker(F_0|_{X_0}\rTo{\ov{\de}_0} (F_1\ot L)|_{X_0}).$$ 
Then Lemma \ref{2-per-res-lem} gives an exact sequence
\begin{equation}\label{exact-functor-eq}
0\to \FF\to (F_1\ot L)|_{X_0}\to\FF'\to 0.
\end{equation}
Since $(F_1\ot L)|_{X_0}\in\Per(X_0)$, this gives 
an isomorphism $\FF'\to\FF[1]$ in $D_{\Sg}(X_0)$, i.e., we obtain a functorial isomorphism
$$t_{\bar{F}}:\fC(\bar{F}[1])\simeq \fC(\bar{F})[1].$$
Now let $f:\bar{E}\to\bar{F}$ be a closed morphism of \mf s.
To see that the image of the triangle \eqref{cone-triangle} under $\fC$ is an exact triangle
in $D_{\Sg}(X_0)$ let us consider the short exact sequence of two-term complexes
\begin{equation}\label{E-F-cone-eq}
0\to [F_1\rTo{\de_1} F_0]\to [C(f)_1\rTo{\de_1} C(f)_0]\to [E[1]_1\rTo{\de_1} E[1]_0]\to 0.
\end{equation}
Since the multiplication by $W$ is injective,
the differentials $\de_1$ in \eqref{E-F-cone-eq} are injective, and so the sequence
of cokernels is exact. 
Thus, we obtain a morphism of short exact sequences of coherent sheaves
\be\label{exact-functor-diagram-eq}
\begin{diagram}
0\rTo{}&\fC(\bar{F})&\rTo{}&\fC(C(f))&\rTo{}&\fC(\bar{E}[1])&\rTo{} 0\\
&\dTo{\id}&&\dTo{\varphi}&&\dTo{\fC(f[1])}\\
0\rTo{}&\fC(\bar{F})&\rTo{}&(F_1\ot L)|_{X_0}&\rTo{}&\fC(\bar{F}[1])&\rTo{} 0
\end{diagram}
\end{equation}
where the second row is the sequence \eqref{exact-functor-eq}, and
$\varphi$ is induced by the morphism
$$F_0\oplus E[1]_0\rTo{(\de,f[1])} F[1]_0|_{X_0}=(F_1\ot L)|_{X_0}.$$
We can extend the first row of \eqref{exact-functor-diagram-eq} to an exact triangle in $D^b(X_0)$
$$\fC(\bar{F})\to \fC(C(f))\to \fC(\bar{E}[1])\rTo{\gamma}\fC(\bar{F})[1].$$
The commutativity of diagram \eqref{exact-functor-diagram-eq}
shows that $\gamma=t_{\bar{F}}\circ\fC(f[1])$,
so the image of the triangle \eqref{cone-triangle} under $\fC$ is an exact 
triangle in $D_{\Sg}(X_0)$. This shows that the functor $\fC$ is exact.
\ed

In the case when $X$ is a smooth affine scheme and $L$ is trivial, the functor $\fC$ is an equivalence
by \cite[Thm.\ 3.9]{Orlov}.
In the non-affine case we need to localize the category $\HMF(X,W)$. Namely, we consider the 
full subcategory
$$\LHZ(X,W)\sub\HMF(X,W)$$
consisting of 
\mf s $\bar{E}$ that are locally contractible
(i.e., there exists an open covering $U_i$ of $X$ in smooth topology such that
$\bar{E}|_{U_i}=0$ in $\HMF(U_i,W|_{U_i})$).
It is easy to see that $\LHZ(X,W)$ is a triangulated subcategory (we will see later that it is thick). 

\begin{defi} For a stack $X$ we define the {\it derived category
of \mf s} by 
$$\DMF(X,W):=\HMF(X,W)/\LHZ(X,W).$$
\end{defi}

The following theorem is a generalization of \cite[Thm.\ 3.9]{Orlov}.

\begin{thm}\label{equiv-thm} 
Let $X$ be a smooth FCDRP-stack (see \ref{stack-def}(iii)),
$L$ a line bundle on $X$ and
$W\in H^0(X,L)$ a potential.
Assume that $W$ is not a zero divisor.
Then the functor 
\begin{equation}\label{main-equiv-eq}
\ov{\fC}:\DMF(W)\to D_{\Sg}(X_0)
\end{equation}
induced by $\fC$ is an equivalence of triangulated categories.
\end{thm}

We will prove this theorem using the following general criterion of equivalence
for triangulated categories.

\begin{lem}\label{equiv-cat-lem} Let $\Phi:\CC\to\DD$ be an exact functor between triangulated
categories. Assume that for every morphism $f:D\to \Phi(C)$, where $C\in\CC$, $D\in\DD$, there
exists an object $C'\in\CC$, a morphism $g:C'\to C$ in $\CC$
and an isomorphism $\a:\Phi(C')\to D$ in $\DD$, such that $\Phi(g)=f\circ\a$.
Then $\Phi$ induces an equivalence $\ov{\Phi}:\CC/\ker(\Phi)\to\DD$ of triangulated categories.
\end{lem}

\Pf . Set $\ov{\CC}=\CC/\ker(\Phi)$. First, let us show that for every $C_1,C_2\in\CC$ the natural
morphism
$$\ov{\Phi}:\Hom_{\ov{\CC}}(C_1,C_2)\to\Hom_{\DD}(\Phi(C_1),\Phi(C_2))$$
is surjective. For $f\in\Hom_{\DD}(\Phi(C_1),\Phi(C_2))$
consider the morphism 
$$\Phi(C_1)\rTo{(\id,f)} \Phi(C_1)\oplus \Phi(C_2)=\Phi(C_1\oplus C_2).$$
By assumption, there exists an object $C'\in\CC$, a morphism
$(q,g):C'\to C_1\oplus C_2$ and an isomorphism $\a:\Phi(C')\to \Phi(C_1)$ such that
$\Phi(q,g)=(\id,f)\circ\a$. In other words, 
we have $\a=\Phi(q)$ and $\Phi(g)=f\circ\a$. Hence, $q$ is an isomorphism in
$\ov{\CC}$ and $\ov{\Phi}(g\circ q^{-1})=f$.
Next, we observe that $\ker(\ov{\Phi})=0$ by definition. 
Now Lemma \ref{faithful-lem} below applied to the functor $\ov{\Phi}$
implies that this functor is faithful.
\ed

\begin{lem}\label{faithful-lem}
Let $\Phi:\CC\to\DD$ be an exact functor between triangulated categories such that $\ker(\Phi)=0$,
and let $C\in\CC$ be a fixed object. Assume that for each $C'\in\CC$
the morphism
$$\Phi:\Hom_{\CC}(C,C')\to\Hom_{\DD}(\Phi(C),\Phi(C'))$$
is surjective. Then this morphism is an isomorphism for every $C'\in\CC$ .
\end{lem}

\Pf .
The argument below is similar to the first part of the proof of \cite[Thm.\ 3.9]{Orlov}).
Suppose $\Phi(f)=0$ for some morphism $f:C\to C'$. Consider
an exact triangle 
$$K\rTo{g} C\rTo{f} C'\to K[1]$$
in $\CC$.
Since the image of this triangle under $\Phi$ is still exact and $\Phi(f)=0$, 
there exists a morphism $h:\Phi(C)\to \Phi(K)$ such that $\Phi(g)\circ h=\id_{\Phi(C)}$.
By assumption we can choose $h'\in\Hom_{\CC}(C,K)$ such that $\Phi(h')=h$. Then
$\Phi(g\circ h')=\id_{\Phi(C)}$, hence, $g\circ h'$ is an isomorphism in $\CC$ (since
$\ker(\Phi)=0$ and $\Phi$ is exact). But $f\circ g\circ h'=0$, so we deduce that $f=0$.
\ed

\begin{lem}\label{loc-free-lem} Let $\FF$ be a MCM sheaf on $X_0$ and let
$p:V\to \FF$ be a surjection of coherent sheaves on $X$, where $V$ is a vector bundle on $X$
and $\FF$ is viewed as a coherent sheaf on $X$.
Then $\ker(p)$ is a vector bundle on $X$.
\end{lem}

\Pf . By taking the pull-back to a presentation of the stack $X$, we can reduce the problem
to the similar question in the case when $X$ is a smooth scheme.
Now we can argue as in the proof of Theorem 3.9 of \cite{Orlov}. Pick a closed point $x\in X$.
Since $X$ is smooth, there exists an integer $N$ such that $\Ext^i_X(\GG,\OO_x)=0$ for $i>N$
for any coherent sheaf $\GG$ on $X$. On the other hand, if $P$ is a vector bundle on $X_0$,
viewed as a coherent sheaf on $X$,
then $\Ext^i_X(P,\OO_x)=0$ for $i>1$. Since by Lemma \ref{MCM-lem}(i),
$\FF$ has a right resolution by vector bundles on $X_0$, we deduce that
$\Ext^i_X(\FF,\OO_x)=0$ for $i>1$. Now the exact sequence
$$0\to \ker(p)\to V\to\FF\to 0$$
implies that $\Ext^i_X(\ker(p),\OO_x)=0$ for $i>0$. It follows that $\ker(p)$ is locally free.
\ed

\noindent
{\it Proof of Theorem \ref{equiv-thm}}.
By Lemma \ref{exact-lem}, $\fC$ is an exact functor. 
To prove that $\ov{\fC}$ is an equivalence
we will apply the criterion of Lemma \ref{equiv-cat-lem} to the functor
$\fC:\HMF(W)\to D_{\Sg}(X_0)$.
Let 
$$f:\FF\to \GG=\fC(\bar{E})$$ 
be a morphism in $D_{\Sg}(X_0)$. We have to construct a \mf\ 
$\bar{F}$, a closed morphism of \mf s $g:\bar{F}\to \bar{E}$ and an isomorphism
$\a:\fC(\bar{F})\to\FF$ such that $\fC(g)=f\circ\a$.
Note that $X_0$ is an FCDRP-stack and has 
finite Krull dimension, as a closed substack of $X$. Also, $X_0$ is Gorenstein as a divisor
in a smooth stack.
Thus, by Lemma \ref{MCM-lem}(iii), we can assume that $\FF$ is a MCM sheaf.
Denote by $E'_i=(E(L^{1/2})_i\ot L^{i/2})|_{X_0}$ the terms of the complex
\eqref{mf-Z-gr-com-eq}. By Lemma \ref{2-per-res-lem}, this complex is exact, so for
large enough $n$ we can apply Lemma \ref{MCM-lem}(iv) to
the exact sequence of sheaves on $X_0$
$$0\to \GG'\to E'_{-n+1}\to\ldots\to E'_{-1}\to E'_0\to\GG\to 0$$
and the morphism $f$.
We obtain an exact sequence of sheaves on $X_0$
$$0\to \FF'\to V_{-n+1}\to\ldots\to V_{-1}\to V_0\to\FF\to 0,$$
and a morphism $f':\FF'\to\GG'$ of sheaves such that
all $V_i$ are vector bundles on $X_0$ and $f'[n]$ represents $f$ in $D_{\Sg}(X_0)$.
Let us prove that from this data one can construct a \mf\ $\bar{F}$, a surjective morphism
$\a:\fC(\bar{F})\to\FF$ of sheaves on $X_0$ and a morphism of \mf s $g:\bar{F}\to\bar{E}$,
such that $\ker(\a)\in\Per(X_0)$ and
$\fC(g)=f\circ\a$ in $\Coh(X_0)$. We use induction in $n$.

If $n=0$ then $f$ is given by a morphism of coherent sheaves $\FF\to\GG$.
We can choose a surjective morphism $p:F_0\to\FF$, where $F_0$ is a vector bundle on $X$,
such that the morphism $f$ extends to a commutative diagram in $\Coh(X)$
\begin{diagram}
F_0 &\rTo{g_0}& E_0\\
\dTo{p}&&\dTo{}\\
\FF&\rTo{f}&\GG
\end{diagram}
By Lemma \ref{loc-free-lem}, $F_1:=\ker(p)$ is a vector bundle. Let $\de_1:F_1\to F_0$
be the natural inclusion.
Since $W\FF=0$, the injective morphism $W:F_0\ot L^{-1}\to F_0$ factors through $\de_1$, so
we obtain an injective morphism $\de_0:F_0\ot L^{-1}\to F_1$ such that $\de_1\de_0=W$.
Note that $\de_1$ and $\de_0$ are isomorphisms over the dense open set $W\neq 0$.
Hence, the equality
$\de_0\de_1\de_0=W\de_0$ implies that $\de_0\de_1=W$, so $\bar{F}=(F_{\bullet},\de_\bullet)$
is a \mf\ of $W$. Furthermore, from the above commutative diagram we get
a unique morphism $g_1:F_1\to E_1$ such that $\de_1g_1=g_0\de_1$. 
It follows that 
$$\de_1\de_0g_0=Wg_0=g_0\de_1\de_0=\de_1 g_1\de_0.$$
Hence, $\de_0g_0=g_1\de_0$, i.e., $g:\bar{F}\to\bar{E}$ is a morphism of \mf s.

Now suppose the assertion is true for $n-1$. Set 
$$\wt{\FF}=\ker(V_0\to \FF) \ \text{ and } 
\wt{\GG}=\ker(E'_0\to\GG).$$
By Lemma \ref{2-per-res-lem},
$\wt{\GG}=\coker((E_0\ot L^{-1})|_{X_0}\to E_1|_{X_0})=\fC(\bar{E}[1])$.
So
we have isomorphisms $\FF\to\wt{\FF}[1]$ and $\GG\to\wt{\GG}[1]$ in $D_{\Sg}(X_0)$, hence
$f$ corresponds to a morphism $\wt{f}:\wt{\FF}\to\wt{\GG}$ in $D_{\Sg}(X_0)$.
By induction assumption, we can construct a \mf\ $\bar{P}$, a surjective morphism
$\a:\fC(\bar{P})\to\wt{\FF}$ and a morphism of \mf s
$\wt{g}:\bar{P}\to\bar{E}[1]$ such that  $\ker(\a)\in\Per(X_0)$
and $\fC(\wt{g})=\wt{f}\circ\a$ in $\Coh(X_0)$.
We have a commutative diagram of sheaves on $X_0$
with exact rows and columns
\begin{diagram}
&&0&&0&&0\\
&&\dTo{}&&\dTo{}&&\dTo{}\\
0&\rTo{}&\ker(\a)&\rTo{\psi}& P_1|_{X_0}&\rTo{}& \KK &\rTo{}&0\\
&&\dTo{}&&\dTo{}&&\dTo{}\\
0&\rTo{}&\fC(\bar{P})&\rTo{\varphi}& P_1|_{X_0}\oplus V_0&\rTo{}&\FF'&\rTo{}&0 \\
&&\dTo{\a}&&\dTo{}&&\dTo{}\\
0&\rTo{}&\wt{\FF}&\rTo{}&V_0&\rTo{}&\FF&\rTo{}&0\\
&&\dTo{}&&\dTo{}&&\dTo{}\\
&&0&&0&&0
\end{diagram}
where $\varphi$ is induced by the natural embedding $\fC(\bar{P})\hra P_1|_{X_0}$
and by the composition $\fC(\bar{P})\rTo{\a}\wt{\FF}\to V_0$, the map $\psi$ is the restriction
of $\varphi$ to $\ker(\a)$, and the sheaves $\KK$ and $\FF'$ are the cokernels of $\psi$
and $\varphi$.
The first row shows that $\KK\in\Per(X_0)$. On the other hand, we have a morphism of exact 
sequences
\begin{diagram}
0&\rTo{}&\fC(\bar{P})&\rTo{}& P_1|_{X_0}\oplus V_0&\rTo{}&\FF'&\rTo{}&0 \\
&&\dTo{}&&\dTo{}&&\dTo{f'}\\
0&\rTo{}&\fC(\bar{P})&\rTo{}& P_1|_{X_0}&\rTo{}&\fC(\bar{P}[-1])&\rTo{}&0
\end{diagram}
Now arguing as in the case $n=0$, we can construct a \mf\ $\bar{F}$, an isomorphism
of sheaves $\fC(\bar{F})\simeq\FF'$ and a closed morphism of \mf s
$\bar{F}\to\bar{P}[-1]$ inducing $f'$. After composing it with $\wt{g}[-1]:\bar{P}[-1]\to\bar{E}$
we get the desired morphism $g:\bar{F}\to \bar{E}$.
Thus, by Lemma \ref{equiv-cat-lem}, we obtain an equivalence of triangulated categories
$$\ov{\fC}:\HMF(W)/\ker(\fC)\to D_{\Sg}(X_0).$$ 

It remains to establish the equality of full subcategories 
$$\ker(\fC)=\LHZ(W)$$
in $\HMF(W)$. By Lemma \ref{quasi-isom-lem}, we have $\LHZ(W)\sub\ker(\fC)$.
Conversely, suppose $\fC(\bar{E})=0$.
Consider an open covering $U_i$ of $X$ in smooth topology, such that $U_i$ are smooth affine
schemes and the pull-backs $L|_{U_i}$ are trivial. By Theorem 3.9 of \cite{Orlov}, the functor $\fC$ induces equivalences 
$$\HMF(U_i,W|_{U_i})\simeq D_{\Sg}(U_i\times_X X_0).$$
Thus, we obtain that $\bar{E}|_{U_i}=0$ in $\HMF(U_i,W|_{U_i})$, and
so $\bar{E}\in \LHZ(W)$.
\ed
 
\begin{cor}\label{mf-nonsing-cor}
Let $X$ be a smooth stack, $W\in H^0(X,L)$ a non-zero-divisor.
The restriction of any \mf\ $\bar{E}\in\HMF(X,W)$ 
to the complement $X\setminus\Sing(X_0)$ of the singular locus of $X_0$
is locally contractible.
\end{cor}

\Pf . We can assume $X$ to be an affine scheme.
Note that the equivalence of Theorem \ref{equiv-thm} is compatible with restrictions to open substacks.
Hence, the assertion follows from the fact that $D_{\Sg}(X\setminus\Sing(X_0))=0$.
\ed

There is one important case when the derived category $\DMF(X,W)$ coincides with the homotopy
category $\HMF(X,W)$.

\begin{prop}\label{affine-case-prop} 
Let us keep the assumptions of Theorem \ref{equiv-thm}.
Assume in addition that $X$ has cohomological dimension $0$. Then
$\ker(\fC)=0$, so in this case 
$$\HMF(X,W)\simeq \DMF(X,W)\simeq D_{\Sg}(X_0).$$
\end{prop}

\Pf . Note that $X_0$ also has cohomological dimension $0$, so every vector bundle on $X_0$ is
a projective object in the category of coherent sheaves on $X_0$.
Thus, we can prove that $\ker(\fC)=0$ by repeating the argument of
\cite[Lem.\ 3.8]{Orlov}. The rest follows from Theorem \ref{equiv-thm}.
\ed

Note that in the case when $X$ is the quotient of an affine scheme by a finite group and $L$ is
trivial, the above result reduces to \cite[Thm 7.3]{Quintero}.
A more general example is the case of a quotient stack $X=U/\Ga$, where $U$ is an affine scheme
and $\Ga$ is a reductive group (see Lemma \ref{PF-lem}(ii)).

\begin{rem}\label{Orlov-def}
Another characterization of $\ker(\fC)$
is given by the unpublished result of Orlov (see \cite[Sec.\ 3.2]{KKP}) that states that $\ker(\fC)$ consists
of all matrix factorizations that appear as direct summands in the convolutions of finite exact
sequences of matrix factorizations. The corresponding quotient category is called 
in \cite{Posic} the {\it absolute derived category}.
\end{rem}

Sometimes it is convenient to use the following generalization of the notion of \mf\ obtained
by replacing vector bundles by coherent sheaves.

\begin{defi}\label{coh-mf-rem} 
A {\it coherent \mf } of a potential $W\in H^0(X,L)$ on a stack $X$ is a pair
of coherent sheaves $F_0$, $F_1$ with maps $\de_1:F_1\to F_0$ and $\de_0:F_0\to F_1$,
such that $\de_0\de_1=W\cdot\id$ and $\de_1\de_0=W\cdot\id$, and
the multiplication by $W$ is injective on both $F_0$ and $F_1$. 
The corresponding homotopy category
$\HMF^c(X,W)$ still has a triangulated structure, and we still have the cokernel functor 
(see \eqref{coker-functor-eq})
$$\fC^c:\HMF^c(X,W)\to D_{\Sg}(X_0)$$ 
which is an exact functor. In the situation of Theorem \ref{equiv-thm}
we can view $\fC^c$ as an exact functor from $\HMF^c(X,W)$ to $\DMF(X,W)$.
\end{defi}

\section{Quasi-\mf s}
\label{quasi-sec}

Here we establish a connection between the category of quasi-\mf s and a
quasicoherent analog of the singularity category. The results of this section will be
used in the study of the push-forward functors for \mf s (see Section \ref{push-forward-sec}).

Assume that $W$ is not a zero divisor.  
Note that we can define the functor 
$$\fC^\infty:\HMF^{\infty}(X,W)\to D'_{\Sg}(X_0)$$ 
for quasi-\mf s in the same way as for \mf s (by taking the cokernel of $\de_1$),
where $D'_{\Sg}(X_0)$ is given by \eqref{quasi-sing-eq}.
It is easy to see that the proof of Lemma \ref{exact-lem} works in this situation,
so the functor $\fC^\infty$ is exact.
Thus, we have a commutative diagram of exact functors 
\begin{diagram}
\HMF(X,W) &\rTo{\fC}&  D_{\Sg}(X_0)\\
\dTo{}&&\dTo{}\\
\HMF^\infty(X,W)&\rTo{\fC^\infty}& D'_{\Sg}(X_0)
\end{diagram}
We are going to prove certain weaker versions of Theorem \ref{equiv-thm} and Proposition
\ref{affine-case-prop} for quasi-\mf s (see Theorem \ref{quasi-mf-thm} below). 

We will need the following fundamental fact about (not necessarily finitely generated)
locally free sheaves.

\begin{lem}\label{proj-mod-lem}
Let $X=\Spec(A)$ be a Noetherian affine scheme, $\FF=\wt{M}$ a quasicoherent
sheaf on $X$ associated with an $A$-module $M$. 
Then $\FF$ is locally free if and only if $M$ is projective.
\end{lem}

\Pf . If $\FF$ is locally free then $M$ is locally projective, so it is projective by \cite[part 2, Sec.\ 3.1]{GR}.
Conversely, assume that $M$ is projective. 
If $M$ is finitely generated then the assertion is
well-known. Otherwise, we can assume that $X$ is connected and apply the result
of Bass \cite[Cor.\ 4.5]{Bass} saying that an infinitely generated projective module over $A$
is free.
\ed

\begin{thm}\label{quasi-mf-thm} 
Let $X$ be a smooth stack with the resolution property, and let $W\in H^0(X,L)$ be a
non-zero-divisor.

(i) Assume that $X=U/\Ga$, where $U$ is an affine scheme and $\Ga$ is a reductive linear
algebraic group. Then the functor $\fC^\infty$ is fully faithful.

(ii) For an arbitrary $X$, let $\LHZ^{\infty}(X,W)\sub \HMF^{\infty}(X,W)$ be the full subcategory  
consisting of quasi-\mf s that are locally homotopic to zero.
Then
$$\ker(\fC^\infty)=\LHZ^\infty(X,W).$$
Moreover, the restriction of an object in $\LHZ^\infty(X,W)$ to any open (in flat topology) 
of the type considered in part (i) is homotopic to zero.
Thus, if we set
$$\DMF^{\infty}(X,W):=\HMF^{\infty}(X,W)/\LHZ^{\infty}(X,W)$$
then the functor $\DMF^\infty(X,W)\to D'_{\Sg}(X_0)$ has zero kernel.

(iii) 
Assume in addition that $X$ has finite cohomological dimension (so $X$ is a smooth FCDRP-stack). Then for $\bar{F}\in\HMF(X,W)$ and $\bar{E}\in\HMF^{\infty}(X,W)$ the map
\begin{equation}\label{mf-inf-map}
\Hom_{\DMF^{\infty}(X,W)}(\bar{F},\bar{E})\to\Hom_{D'_{\Sg}(X_0)}(\fC(\bar{F}),\fC^\infty(\bar{E}))
\end{equation}
is an isomorphism.
\end{thm}

For the proof we need the following analogs of Lemma \ref{MCM-lem}(ii) and \cite[Prop.\ 1.21]{Orlov}.

\begin{lem}\label{qc-loc-free-lem} 
Let $A$ be a Noetherian Gorenstein commutative ring of Krull dimension $n$. 

(i) For any projective $A$-module $P$ and an $A$-module $M$ one has
$\Ext_A^i(M,P)=0$ for $i>n$.

(ii) Let $\FF=\wt{M}$ be a quasicoherent sheaf on $X=\Spec(A)$
associated with an $A$-module $M$, that admits a right
locally free resolution $\FF\to Q^\bullet$ and such that $\FF\in\Lfr(X)$.
Then $M$ is a projective $A$-module and so, by Lemma \ref{proj-mod-lem},
$\FF$ is a locally free sheaf.
\end{lem} 

\Pf . (i) It is enough to consider the case when $P$ is a free module.
Since $A$ is Gorenstein of Krull dimension $n$, the injective dimension of $A$ is equal to $n$.
Since (infinite) direct sums of injective modules over a Noetherian ring are also injective,
it follows that the injective dimension of $P$ is also $n$.

(ii) By assumption there exists a finite locally free resolution
$$0\to Q^{-N}\to\ldots\to Q^{-1}\to\FF\to 0.$$
Sewing it with the right resolution $Q^0\to Q^1\to\ldots $ of $\FF$ we obtain
the exact complex
$$0\to Q^{-N}\to\ldots\to Q^{-1}\to Q^0\to Q^1\to\ldots$$
of locally free sheaves. By Lemma \ref{proj-mod-lem}, 
the $A$-modules corresponding to $Q^i$ are projective.
Thus, it suffices to prove that for any infinite exact complex
$$0\to P^0\to P^1\to\ldots$$
of projective modules, the module $N=\coker(P^0\to P^1)$ is also projective.
By part (i), any projective $A$-module also has
the injective dimension $\le n$.
Let us set $N^i=\coker(P^{i-1}\to P^i)$, where $i\ge 1$, so that $N=N^1$.
Then using the exact sequences $0\to N^i\to P^{i+1}\to N^{i+1}\to 0$ we obtain
$$\Ext^1(N,P^0)=\Ext^2(N^2,P^0)=\ldots=\Ext^{n+1}(N^{n+1},P^0)=0.$$ 
Hence, the sequence
$$0\to P^0\to P^1\to N\to 0$$
splits, so $N$ is projective. 
\ed

\begin{lem}\label{qc-mf-sur-lem} 
Let $X=U/\Ga$ be a nice affine Gorenstein quotient stack, where $U$ is a Noetherian affine scheme
of finite Krull dimension $n$, and $\Ga$ is a reductive linear algebraic group. 

(i) If $\FF$ is a quasicoherent sheaf and $\PP$ is a locally free sheaf on $X$ then
$\Ext^i_X(\FF,\PP)=0$ for $i>n$.

(ii) Let $\FF$ and $\GG$ 
be quasicoherent sheaves on $X$ such that $\FF$ admits a right locally free resolution 
$\FF\to Q^\bullet$. Then the natural map
$$\Hom_X(\FF,\GG)\to\Hom_{D'_{\Sg}(X)}(\FF,\GG)$$
is surjective.
\end{lem}

\Pf . (i) By Lemma \ref{PF-lem}(i), it is enough to prove a similar statement on $U$. But $U=\Spec(A)$ is affine and Noetherian, so $\PP$ corresponds to a projective $A$-module by Lemma \ref{proj-mod-lem}.
Therefore, it has injective dimension $\le n$ by Lemma \ref{qc-loc-free-lem}(i).

(ii) The proof is analogous to that of \cite[Prop.\ 1.21]{Orlov}, using part (i) together with
the fact that for any locally
free sheaf $\PP$ and any quasicoherent sheaf $\FF$ on $X$ one has
$\Ext^i_X(\PP,\FF)=0$ for $i>0$. Indeed, by Lemma 
\ref{PF-lem}(i), this reduces to a similar statement on $U$, which follows from Lemma \ref{proj-mod-lem}.
\ed

\noindent
{\it Proof of Theorem \ref{quasi-mf-thm}.}
(i) First, as in the proof of Theorem \ref{equiv-thm}, we see that the functor $\fC^\infty$ is exact.
Next, we claim that in this case $\ker(\fC^\infty)=0$.
Indeed, suppose $\bar{P}\in\MF^\infty(X,W)$ is such that the quasicoherent
sheaf $\FF=\coker(P_1\to P_0)$ belongs to $\Lfr(X_0)\sub D^b(\Qcoh(X_0))$. 
By Lemma \ref{2-per-res-lem}, the sheaf $\FF$ has a right locally free resolution.
Hence, we can apply Lemma \ref{qc-loc-free-lem}(ii) to an affine covering of $X$ 
to deduce that $\FF$ is locally free. 
Note that $X_0$ itself is a quotient of a Noetherian affine scheme by $\Ga$.
In particular, $X_0$ has cohomological dimension $0$, so by Lemma \ref{PF-lem}(ii), 
$\FF$ is a projective object in $\Qcoh(X_0)$. Therefore, 
the same argument as in \cite[Lem.\ 3.8]{Orlov}, shows that $\bar{P}$ is contractible.
Now Lemma \ref{qc-mf-sur-lem}(ii) and the projectivity of locally free sheaves on $X$ imply that
the functor $\fC^\infty$ is full (using an argument of \cite[Lem.\ 3.5]{Orlov}).
By Lemma \ref{faithful-lem}, the functor $\fC^\infty$ is also faithful.

(ii) This follows from (i). 

(iii) By Lemma \ref{faithful-lem}, it is enough to prove that the map \eqref{mf-inf-map}
is surjective.
To this end we use the strategy similar to that of the proof of Theorem \ref{equiv-thm}.
Let us set $\FF=\fC(\bar{F})$ (this is a coherent sheaf on $X_0$), 
and suppose we have a morphism
$f:\FF\to\fC^\infty(\bar{E})$ in $D'_{\Sg}(X_0)$. 
Let us consider the quasicoherent sheaf on $X_0$
$$\GG=\FF\oplus\fC^\infty(\bar{E})=\fC^\infty(\bar{F}\oplus\bar{E})$$ 
and the morphism $(\id,f):\FF\to \GG$.
Applying the same inductive procedure as in the proof of Theorem \ref{equiv-thm}
(using Lemma \ref{MCM-lem}(iv)) to the morphism $(\id,f)$ we can construct a \mf\ 
$\bar{F'}\in\MF(X,W)$, a closed morphism of quasi-\mf s $g':\bar{F'}\to\bar{F}\oplus\bar{E}$ and an isomorphism $\a:\fC(\bar{F'})\to\FF$ in $D_{\Sg}(X_0)$ such that $\fC^\infty(g')=(\id,f)\circ\a=(\a,f\circ\a)$.
Let us write $g'=(q,g)$, where $q:\bar{F'}\to\bar{F}$ and $g:\bar{F'}\to\bar{E}$.
Then $\fC(q)=\a$ and $\fC^\infty(g)=f\circ \a$.
Since $\a$ is an isomorphism, from part (ii) we obtain that $q$ is an isomorphism in $\DMF^\infty(X,W)$.
Hence, the morphism $g\circ q^{-1}:\bar{F}\to\bar{E}$ in $\DMF^\infty(X,W)$ is in the preimage
of $f$ under $\fC^\infty$.
\ed

\begin{cor} In the situation of Theorem \ref{quasi-mf-thm}(i) the natural functor
$$\HMF(X,W)\to\HMF^{\infty}(X,W)$$
is fully faithful.
\end{cor}

\begin{lem}\label{idemp-lem} The triangulated category $\DMF^\infty(X,W)$ has infinite direct sums. 
Hence, all idempotents in this category are split.
\end{lem}

\Pf . Lemma \ref{proj-mod-lem} implies that an (infinite) direct sum of locally free sheaves on $X$
is still locally free. Now, the standard construction for the homotopy category of complexes
can be adapted to show that infinite direct sums exist in $\HMF^{\infty}(X,W)$. 
By Theorem \ref{quasi-mf-thm},
the subcategory $\LHZ^\infty=\ker(\fC^\infty)$ is closed under arbitrary direct sums.
By \cite[Lem.\ 1.5]{BN}, this implies that the quotient category 
$\DMF^\infty(X,W)$ has direct sums. The last assertion follows by \cite[Prop.\ 3.2]{BN}.
\ed

Abusing  notation we will also denote by 
$$\fC^\infty:\DMF^\infty(X,W)\to D'_{\Sg}(X_0)$$
the exact functor induced by $\fC^\infty$.

By Lemma \ref{idemp-lem}, the natural functor $\DMF(X,W)\to\DMF^\infty(X,W)$ extends to
the functor
$$\iota:\ov{\DMF}(X,W)\to\DMF^\infty(X,W).$$

\begin{prop}\label{qc-isom-prop} 
Assume that $X$ is a smooth FCDRP-stack,
and let $W\in H^0(X,L)$ be a non-zero-divisor.
Then for
objects $(\bar{E},e)\in\ov{\DMF}(X,W)$ and $\bar{F}\in\DMF^{\infty}(X,W)$ 
the morphism
$$\fC^\infty:\Hom_{\DMF^{\infty}(X,W)}(\iota(\bar{E},e),\bar{F})\to
\Hom_{D'_{\Sg}(X_0)}(\fC^\infty(\iota(\bar{E},e)),\fC^\infty(\bar{F}))$$
is an isomorphism. Furthermore, if $\a\in\Hom_{\DMF^{\infty}(X,W)}(\iota(\bar{E},e),\bar{F})$
is such that $\fC^\infty(\a)$ is an isomorphism in $D'_{\Sg}(X_0)$
then $\a$ is an isomorphism in $D'_{\Sg}(X_0)$.
\end{prop}

\Pf . It is enough to check the first assertion for $\bar{E}$ instead of $\iota(\bar{E},e)$, in which
case it follows from Theorem \ref{quasi-mf-thm}(iii).
The second assertion follows from the fact that the exact functor
$\fC^\infty:\DMF^\infty(X,W)\to D'_{\Sg}(X_0)$ has zero kernel by Theorem \ref{quasi-mf-thm}(ii).
\ed

Quasi-\mf s provide a natural setup for defining the functors of push-forwards with respect to
smooth affine morphisms with geometrically integral fibers. Recall that by \cite[part 1, (3.3.1)]{GR},
the push-forward of a locally free sheaf under such a morphism  is locally projective, hence,
by Lemma \ref{proj-mod-lem}, locally free. This allows us to make the following definition.

\begin{defi}\label{quasi-mf-push-forward} 
Let $f:X\to Y$ be a smooth affine morphism of stacks with geometrically integral fibers,
let $W\in H^0(Y,L)$ be a potential, and
let $\bar{E}=(E,\de)$ be a quasi-\mf\  of $f^*W$.
The {\it push forward quasi-\mf\ } $f_*\bar{E}$ of $W$ is 
the pair of locally free sheaves $(f_*E_0, f_*E_1)$ together with the differential $f_*\de$.
\end{defi}

This gives a dg-functor $f_*:\MF^\infty(X,f^*W)\to \MF^\infty(Y,W)$, which induces exact functors
$$f_*:\HMF^\infty(X,f^*W)\to\HMF^\infty(Y,W).$$

Assume in addition that $W$ is not a zero divisor, and let $Y_0$ (resp., $X_0$) be the zero locus
of $W$ (resp., of $f^*W$). Let $g:X_0\to Y_0$ denote the morphism induced by $f$.
Since the morphism $X_0\to Y_0$ is also smooth with geometrically integral fibers,
the push-forward functor $g_*$ on quasicoherent sheaves induces
an exact functor $D^b(\Qcoh(X_0))\to D^b(\Qcoh(Y_0))$ that sends $\Lfr(X_0)$ to $\Lfr(Y_0)$,
and hence gives an exact functor $g_*:D'_{\Sg}(X_0)\to D'_{\Sg}(Y_0)$.
Furthermore, since $f$ is an affine morphism, for every $\bar{E}\in \MF^\infty(X,f^*W)$
we have a natural isomorphism of quasicoherent sheaves on $X_0$
$$g_*\coker(E_1\rTo{\de} E_0)\simeq \coker(f_*E_1\rTo{f_*\de} f_*E_0).$$
Thus, we obtain a commutative diagram
of exact functors between triangulated categories
\begin{diagram}
\HMF^\infty(X,f^*W)&\rTo{\fC^\infty}&D'_{\Sg}(X_0)\\
\dTo{f_*}&&\dTo{g_*}\\
\HMF^\infty(Y,W)&\rTo{\fC^\infty}&D'_{\Sg}(Y_0)
\end{diagram}
If $Y$ is smooth and has the resolution property, then by Theorem \ref{quasi-mf-thm}
we obtain an induced functor of derived categories
$$f_*:\DMF^\infty(X,f^*W)\to \DMF^\infty(Y,W)$$
also compatible with the functor $g_*:D'_{\Sg}(X_0)\to D'_{\Sg}(Y_0)$.

\section{Supports}
\label{supports-sec}

Orlov showed in \cite{O-compl} that objects of the singularity category $D_{\Sg}(X)$ can
be represented as direct summands of complexes with cohomology supported on the
singular locus of $X$. In this section we will introduce and study a more general notion
of support for objects of the singularity category and the corresponding notion of support
for \mf s.

\begin{defi}
Let $X$ be a Gorenstein stack with the resolution property, and let
$Z\sub X$ be  a closed substack. We define the singularity category of $X$ with support on $Z$ as
$$D_{\Sg}(X,Z)=D^b(X,Z)/\perf(X,Z),$$
where $D^b(X,Z)$ is the triangulated 
subcategory of $D^b(X)$ of complexes with cohomology supported
on $Z$, and $\perf(X,Z)\sub D^b(X,Z)$ consists of perfect complexes in $D^b(X,Z)$.
We denote by $\ov{D_{\Sg}(X,Z)}$ the 
idempotent completion of $D_{\Sg}(X,Z)$.
\end{defi}

By \cite[Lem.\ 2.6]{O-compl}, $D_{\Sg}(X,Z)$ is a full subcategory of $D_{\Sg}(X)$.

The following result is a slight generalization of
\cite[Prop.\ 2.7]{O-compl} 
and is proved similarly (see also \cite[Thm.\ 1.3]{Chen}).

\begin{prop}\label{direct-summand-prop} Assume that $X$ is a Gorenstein stack of finite Krull
dimension with the resolution property. 
Let $Z\sub X$ be a closed substack and let $j:U=X\setminus Z\to X$ be the open embedding
of its complement.
Then the subcategory $\ov{D_{\Sg}(X,Z)}\sub \ov{D_{\Sg}(X)}$ coincides with
the kernel of the functor
$$\ov{j^*}:\ov{D_{\Sg}(X)}\to\ov{D_{\Sg}(U)},$$
induced by the restriction functor $j^*$.
\end{prop}

\Pf . 
First, let us consider the functor
$$j^*:D_{\Sg}(X)\to D_{\Sg}(U).$$
We claim that $\ker(j^*)$ consists of direct summands of objects in $D_{\Sg}(X,Z)$.
Indeed, it is clear that for $F\in D_{\Sg}(X,Z)$ one has $j^*F=0$.
Conversely, suppose that for $F\in D^b(X)$ we have $j^*F=0$ in $D_{\Sg}(U)$, i.e, $j^*F$
is a perfect complex.
By  Lemma \ref{MCM-lem}(iii), we can assume that $F$ is a MCM-sheaf.
The condition that $j^*F$ is perfect implies by Lemma \ref{MCM-lem}(ii)
that $j^*F$ is a vector bundle. We have to show that $F$ is a direct summand of an object in
$D_{\Sg}(Y)$ represented by a complex with cohomology supported in $Z$.
Note that $U$ has finite cohomological dimension as an open substack of $X$ (see Lemma 
\ref{FCD-stack-lem}).
Take a resolution of $F$ by vector bundles
$$\ldots P_2\to P_1\to P_0\to F\to 0$$
and consider the sheaf $G=\ker(P_n\to P_{n-1})$, where $n$ is the
the cohomological dimension of $U$.
We have the corresponding morphism $\a:F\to G[n+1]$ in $D(X)$. Since
the sheaf $j^*F$ is a vector bundle, and the cohomological dimension of $U$ is $n$,
we have $j^*\a=0$. Hence, $\a$ factors through an object
$A$ of $D(X,Z)$. But $\a$ descends to an isomorphism in $D_{\Sg}(X)$, so
our claim follows. 

Next, we observe that $D_{\Sg}(U)$ can be identified with the quotient category
$D_{\Sg}(X)/\ker(j^*)$. Indeed, this can be easily deduced from the fact that every object (resp.,
morphism) in $D^b(U)$ can be extended to an object (resp., morphism) in $D^b(X)$ (see
\cite[Lem.\ 2.12]{AB}. Thus, the first part of the argument implies that every morphism
in $D_{\Sg}(X)$ that becomes zero in $D_{\Sg}(U)$ factors through an object of $D_{\Sg}(X,Z)$.

Now suppose we have an object $(F,e)\in\ov{D_{\Sg}(X)}$, where $F\in D_{\Sg}(X)$ and 
$e:F\to F$ is a projector, such that $\ov{j^*}(F,e)=0$. Then the morphism $e:F\to F$ becomes
zero in $D_{\Sg}(U)$, hence it factors through an object $A\in D_{\Sg}(X,Z)$. But
the identity morphism of $(F,e)$ factors through $e$ in $\ov{D_{\Sg}(X)}$, so
$(F,e)$ is a direct summand of $A$ in $\ov{D_{\Sg}(X)}$.
\ed

This proposition implies a version of \cite[Prop.\ 2.7]{O-compl} for stacks.

\begin{cor}\label{sing-supp-cor} Let $X$ be a Gorenstein stack of finite Krull
dimension with the resolution property. 
The natural inclusion 
$$\ov{D_{\Sg}(X,\Sing(X))}\to \ov{D_{\Sg}(X)}$$
is an equivalence of categories.
\end{cor}

Now we are going to introduce the notion of support for \mf s.
Let  $W\in H^0(X,L)$ be a non-zero-divisor potential on a smooth FCDRP-stack $X$,
and 
consider its zero locus $X_0=W^{-1}(0)$.
For a closed substack $Z\sub X_0$
we would like to characterize
the subcategory in $\DMF(W)$ corresponding to $D_{\Sg,Z}(X_0)$ under
the equivalence of Theorem \ref{equiv-thm}.
Given a matrix factorization $\bar{P}=(P,\de)$ of $W$ and a closed point $x$ in a presentation of 
$X_0$, consider the complex 
$$(i_x^*P,\de_x):=i_x^*\com(\bar{P})$$
(see \eqref{mf-Z-gr-com-eq}), where $i_x:x\to X$ is the natural embedding. Denote 
by $H^*(\bar{P},x)$ the cohomology of this complex.

\begin{lem}\label{point-lem} 
Let $\bar{P}=(P,\de)$ be a \mf\ of $W$ on $X$, where $X$ and $W$ are as above.

(i) If in addition $X$ is an affine scheme then for any
coherent sheaf  $\GG$ on $X_0$ 
$$\Hom_{D_{\Sg}(X_0)}(\fC(\bar{P}),\GG)\simeq H^0(\GG\ot \com(\bar{P})^\vee),$$
where $\fC(\bar{P})=\coker(\de:P_1\to P_0)$.

(ii) Assume that $X$ is an affine scheme.
For any closed point $x\in X_0$ one has natural isomorphisms
$$H^i(\bar{P},x)^*\simeq\Hom_{D_{\Sg}(X_0)}(\fC(\bar{P}),\OO_x[i]), \ \ i=0,1.$$ 

(iii) For any closed point $x$ in a presentation of $X_0$ 
one has $H^0(\bar{P},x)=0$ if and only if $H^1(\bar{P},x)=0$ if and only if
$\fC(\bar{P})$ is locally free on $X_0$ in a neighborhood of $x$. Equivalently,
there exists an open neighborhood $U$ of $x$ in $X$ such that
$\bar{P}|_U=0$ in $\HMF(U,W|_U)$, i.e., $\bar{P}|_U$ is a contractible \mf.
\end{lem}

\Pf . (i) Since $X$ is affine, Lemma 3.6 and Proposition 1.21 of \cite{Orlov} imply that 
$$\Hom_{D_{\Sg}(X_0)}(\fC(\bar{P}),\GG)\simeq\Hom_{X_0}(\fC(\bar{P}),\GG)/\RR$$
where $\RR$ consists of morphisms that factor through a vector bundle on $X_0$.
By Lemma \ref{2-per-res-lem}, the complex $\com(\bar{P})$ 
is exact, so we can identify
$\fC(\bar{P})$ with a subsheaf in $(P_1\ot L)|_{X_0}$.
We claim that in fact the subspace $\RR$ coincides with the image of the natural map
$\Hom(P_1\ot L|_{X_0},\GG)\to\Hom(\fC(P),\GG)$.
Indeed, if $V$ is a vector bundle on $X_0$ then a map $\phi:\fC(P)\to V$ defines 
an element in 
$$\psi\in\ker\left(\Hom_{X_0}(P_0|_{X_0},V)\to\Hom_{X_0}(P_1|_{X_0},V)\right).$$
Since $X_0$ is affine and the complex $\com(\bar{P})$ of vector bundles on $X_0$ is exact, 
the complex $\Hom_{X_0}(\com(\bar{P}),V)$ is also exact. Therefore,
$\psi$ comes from an element of
$\Hom_{X_0}(P_1\ot L|_{X_0},V)$, i.e., $\phi$ factors through a map $P_1\ot L|_{X_0}\to V$.
Thus, we get an isomorphism
$$\Hom_{D_{\Sg}(X_0)}(\fC(\bar{P}),\GG)\simeq
\coker\left(\Hom_{X_0}(P_1\ot L|_{X_0},\GG)\to\Hom_{X_0}(\fC(\bar{P}),\GG)\right).$$
It remains to use the isomorphism 
$$\Hom_{X_0}(\fC(\bar{P}),\GG)\simeq\ker\left(
\Hom_{X_0}(P_0|_{X_0},\GG)\to\Hom_{X_0}(P_1|_{X_0},\GG)\right).$$

(ii) 
By part (i), we have an isomorphism
$$H^0(\bar{P},x)^*\simeq\Hom_{D_{\Sg}(X_0)}(\fC(\bar{P}),\OO_x).$$
The statement for $H^1(\bar{P},x)$ follows because $\fC$ is
an exact functor. 

(iii) We can assume that $X$ is an affine scheme.
Recall that $\FF=\fC(\bar{P})$ is locally free if and only it is a zero object of $D_{\Sg}(X_0)$
(by Lemmas 1.20 and 3.6 of \cite{Orlov}). 
Hence, in this case $H^0(P,x)=0$. Conversely,
assume that $H^0(P,x)=0$. Consider the exact sequence
$$0\to \FF \rTo{\a} P_1|_{X_0}\to\FF'\to 0,$$ 
where $\FF'=\coker(\de:P_0\to P_1)$. Then it is easy to see that 
$$H^0(\bar{P},x)=\ker(\a(x): i_x^*\FF\to i_x^*P_1).$$ 
Thus, the vanishing of $H^0(\bar{P},x)$ implies that
$\a(x)$ is injective, i.e, $\Tor^1(\FF',\OO_x)=0$. This implies that
$\FF'$ is locally free, hence, $\FF$ is locally free (since $\FF'\simeq\FF[1]$ in $D_{\Sg}(X_0)$).
Finally, we observe that
the ranks of $H^0(\bar{P},x)$ and of $H^1(\bar{P},x)$ are equal (since the ranks of $P_0$ and $P_1$
are equal).
\ed

\begin{defi} Let $X$ be a stack and $W\in H^0(X,L)$.
For a closed substack $Z\sub X_0=W^{-1}(0)$ 
the category of {\it \mf s of $W$ with support on $Z$}
is the full subcategory $\HMF(X,Z;W)\sub\HMF(X,W)$ consisting of
$\bar{P}\in \HMF(X,Z;W)$ such that
$$H^*(\bar{P},x)=0\ \text{ for all } x\in \wt{X_0\setminus Z}.$$
Here $x$ runs over closed points of some presentation $\wt{X_0\setminus Z}\to X_0\setminus Z$.
Note that $\LHZ(X,W)$ is a subcategory of $\HMF(X,Z;W)$. 

We define the {\it derived category of \mf s with support on }
$Z$ by
$$\DMF(X,Z;W):=\HMF(X,Z;W)/\LHZ(X,W),$$
and denote its idempotent completion by $\ov{\DMF}(X,Z,W)$. 
\end{defi}

Note that $\DMF(X,Z;W)$ (resp., $\ov{\DMF}(X,Z;W)$)
is a full triangulated subcategory in $\DMF(X,W)$ (resp., in $\ov{\DMF}(X,W)$). Furthermore,
the functors 
$$H^*(?,x):\DMF(X,W)\to k-\mod$$ 
extend naturally to the category $\ov{\DMF}(X,W)$.
Thus, we have
\be\label{idemp-H*-eq}
\ov{\DMF}(X,Z;W)=\bigcap_{x\in\wt{X_0\setminus Z}}\ker H^*(?,x).
\end{equation}

\begin{prop}\label{support-prop} 
Let $X$ be a smooth FCDRP-stack. 
The equivalence
$$\ov{\fC}:\ov{\DMF}(X,W)\to\ov{D_{\Sg}(X_0)}$$
induced by the functor \eqref{main-equiv-eq} identifies the full subcategories
$\ov{\DMF}(X,Z;W))\sub \ov{\DMF}(X,W)$ and
$\ov{D_{\Sg}(X_0,Z)}\sub \ov{D_{\Sg}(X_0)}$.
\end{prop}

\Pf . 
By Proposition \ref{direct-summand-prop}, it suffices to check that 
$\ov{\DMF}(X,Z;W)$ coincides with the kernel of the restriction functor
$$\ov{j^*}:\ov{\DMF}(X,W)\to\ov{\DMF}(X\setminus Z,W|_{X\setminus Z}).$$
The fact that the subcategory
$\ov{\DMF}(X,Z;W)$ contains $\ker \ov{j^*}$ follows immediately from \eqref{idemp-H*-eq}. 
To prove the opposite inclusion we have to check that
$j^*(\DMF(X,W))=0$. 
But this follows from Lemma \ref{point-lem}(iii) because vanishing
of an object in $\DMF(X\setminus Z,W|_{X\setminus Z})$ 
is a local property. 
\ed

\section{Push-forwards}
\label{push-forward-sec}

Here we apply the results and constructions of Sections 4 and 5 to
define and study the push-forward functors for categories of \mf s.
We will use these functors in \cite{PV-CohFT} to give an algebraic
construction of a cohomological field theory related to the Landau-Ginzburg
model for a quasihomogeneous isolated singularity (see \cite{FJR}).

\begin{prop}\label{push-forward-prop}
Let $f:X\to Y$ be a representable morphism of
smooth FCDRP-stacks and $W\in H^0(Y,L)$ a potential
such that $W$ and $f^*W$ are not zero divisors.
Let $Z\sub X_0$ be a closed substack of the zero locus of $f^*W$, such that
the induced morphism $f:Z\to Y$ is proper. Let $Y_0\sub Y$ denote the zero locus of $W$, and
let $f_0:X_0\to Y_0$ be the map induced by $f$. 

(i)
The derived push-forward functor
$$Rf_{0*}:D^b(X_0,Z)\to D^b(Y_0,f(Z))$$
induces a functor
$$\ov{D_{\Sg}(X_0,Z)}\to\ov{D_{\Sg}(Y_0,f(Z))},$$
and hence, by Proposition \ref{support-prop}, a functor
\begin{equation}\label{mf-push-forward}
Rf_*:\ov{\DMF}(X,Z;f^*W)\to \ov{\DMF}(Y,f(Z);W).
\end{equation}

(ii) Assume that $f$ is flat. Then
for $\bar{F}\in\ov{\DMF}(X,Z;f^*W)$ and $\bar{E}\in\ov{\DMF}(Y,W)$ one has a natural isomorphism
$$\Hom(\bar{E},Rf_*\bar{F})\simeq\Hom(f^*\bar{E},\bar{F}).$$

(iii)
If we have a representable morphism $f':Y\to S$ to a smooth FCDRP-stack such that
$W$ is a pull-back of a potential $W'$ on $S$ and $f(Z)$ is proper over $S_0$ then
we have an isomorphism of functors
$$R(f'\circ f)_*\simeq Rf'_*\circ Rf_*$$ 
from $\ov{\DMF}(X,Z;f^*W)$ to $\ov{\DMF}(S,W)$.
\end{prop}

\Pf . (i) We have to check that if $V$ is a
perfect complex on $X_0$ with cohomology supported on $Z$ 
then $Rf_{0*}(V)$ is a perfect complex on $Y_0$. By passing to a presentation of $Y_0$ we
can work with schemes.
By \cite[Cor.\ 5.8.1]{SGA6I}, it suffices to check that $Rf_{0*}(V)$ has bounded coherent cohomology
and is of finite tor-dimension.
The former assertion follows from the fact that $Z$ is proper over $Y$, so
by \cite[Cor.\ 3.7.2]{SGA6III}, it is enough to show that the map $g$ is of finite tor-dimension.
Note that the assumption that $f^*W$ is not a zero divisor implies that
  the Cartesian diagram
\begin{diagram}
X_0 &\rTo{} & X\\
\dTo{f_0} && \dTo{f}\\
Y_0 & \rTo{} & Y
\end{diagram}
is tor-independent. Since $f$ is of finite tor-dimension (as a map between smooth stacks), it follows
that $f_0$ is also of finite tor-dimension. 

(ii) Since $f$ is flat, the pull-back functor $f^*$ is compatible with the
equivalences of Theorem \ref{equiv-thm}. Now the statement follows from the
similar isomorphism for derived categories of sheaves (cf. \cite[Lem.\ 1.2]{Orlov}).

(iii) This follows immediately from the similar property of the push-forward functors for derived
categories of sheaves.
\ed

\begin{rem}\label{coh-push-forward-rem}
If in the situation of Proposition \ref{push-forward-prop}(i) the morphism $f$ is proper then we
do not need idempotent completions. Just using the equivalences of Theorem \ref{equiv-thm}
gives a functor
$$Rf_*:\DMF(X,f^*W)\to \DMF(Y;W).$$
This functor sends \mf s with support on $Z$ to \mf s with support on $f(Z)$, as
one can see directly from Lemma \ref{point-lem}(iii).
Furthermore, if $f$ is a finite morphism then the above functor is induced by the composition
$$f_*:\HMF(X,f^*W)\to \HMF^c(Y,W)\rTo{\fC^c} \D_{\Sg}(Y_0)\simeq\DMF(Y,W),$$
where the first arrow sends a \mf\ $(P_\bullet,\de)$ of $f^*W$ to the coherent \mf\ 
$(f_*P_\bullet,f_*\de)\in\HMF^c(Y,W)$ (see Definition \ref{coh-mf-rem}).
\end{rem}

Note that now we have two different
notions of the push-forward functors: one for quasi-\mf s (see Definition
\ref{quasi-mf-push-forward}) 
and another given by the above proposition. The next result shows that they are compatible.

\begin{prop}\label{compat-push-forward-prop} 
Let $X$ and $Y$ be smooth FCDRP-stacks,
$f:X\to Y$ a smooth affine morphism with geometrically integral fibers,
$W\in H^0(Y,L)$ a non-zero-divisor,
and $Z\sub X_0$ a closed substack of the zero locus of $W$.
Then the following diagram of functors is commutative
\begin{equation}\label{quasi-mf-push-forward-diagram}
\begin{diagram}
\ov{\DMF}(X,Z;f^*W) & \rTo{f_*}& \ov{\DMF}(Y,W)\\
\dTo{}&&\dTo{}\\
\DMF^\infty(X,f^*W) &\rTo{f_*}& \DMF^\infty(Y,W)
\end{diagram}
\end{equation}
\end{prop}

\Pf . Let $f_0:X_0\to Y_0$ be the morphism induced by $f$.
We have a commutative diagram of push-forward functors
\begin{equation}\label{sing-push-forward-diagram}
\begin{diagram}
\ov{D_{\Sg}(X_0,Z)} &\rTo{f_{0*}}& \ov{D_{\Sg}(Y_0)}\\
\dTo{}&&\dTo{}\\
D'_{\Sg}(X_0) &\rTo{f_{0*}}& D'_{\Sg}(Y_0)
\end{diagram}
\end{equation}
Furthermore, each of the arrows in the diagram \eqref{quasi-mf-push-forward-diagram}
is compatible with the corresponding arrow in the diagram \eqref{sing-push-forward-diagram}
via the appropriate functors $\fC$ or $\fC^\infty$.
Now the assertion follows from Proposition \ref{qc-isom-prop}. 
\ed

{\sc Department of Mathematics, University of Oregon, Eugene, OR 97405}

{\it Email addresses}: apolish@uoregon.edu, vaintrob@uoregon.edu

\end{document}